\documentclass[11pt,oneside, reqno, twoside]{amsart}
\pdfoutput=1 
\usepackage{accents}
\usepackage[portrait, margin=1in]{geometry}
\usepackage[OT2,T1]{fontenc}
\usepackage{amsmath}
\usepackage{amsthm}
\usepackage{amssymb}
\usepackage{mathtools}
\usepackage{mathrsfs}
\usepackage[shortlabels]{enumitem}
\usepackage{graphicx}
\usepackage{float}

\usepackage{fancyhdr}
\def\zz{\mathbb{Z}}

\def\mcA{\mathcal A}

\def\mcB{\mathcal B}

\def\mcE{\mathcal{E}}

\def\ff{\mathbb{F}}

\theoremstyle{plain}\newtheorem{thm}{Theorem}[section]
\theoremstyle{plain}\newtheorem*{thm*}{Theorem}
\theoremstyle{plain}
\theoremstyle{plain}
\theoremstyle{plain}\newtheorem{conjecture}[thm]{Conjecture}
\theoremstyle{plain}\newtheorem{lemma}[thm]{Lemma}
\theoremstyle{plain}
\theoremstyle{plain}\newtheorem{remark}[thm]{Remark}
\theoremstyle{plain}

\usepackage[style=alphabetic, maxnames=50]{biblatex}

\usepackage{fullpage}
\usepackage{xstring}

\usepackage{hyperref}
\usepackage{dsfont}
\usepackage{indentfirst}

\hypersetup{
    colorlinks=true,
    linkcolor=blue,
    filecolor=magenta,      
    urlcolor=cyan,
    pdftitle={Numerical Inverstigation of Lower Order Biases in Moment Expansions of One Parameter Families of Elliptic Curves},
    pdfpagemode=FullScreen,
    citecolor=blue,
    pdfauthor={Timothy Cheek, Pico Gilman, Kareem Jaber, Steven J. Miller, Marie-H\'el\`ene Tom\'e}
}

\bibliography{references.bib}

\title{Numerical Investigation of Lower Order Biases in Moment Expansions of One Parameter Families of Elliptic Curves}

\author{Timothy Cheek}
\address{Department of Mathematics, University of Michigan}
\email{\href{mailto:timcheek@umich.edu}{timcheek@umich.edu}}

\author{Pico Gilman}
\address{Department of Mathematics, University of California Santa Barbara}
\email{\href{mailto:picogilman@gmail.com}{picogilman@gmail.com}}

\author{Kareem Jaber}
\address{Department of Mathematics, Princeton University}
\email{\href{mailto:kj5388@princeton.edu}{kj5388@princeton.edu}}

\author{Steven J. Miller}
\address{Department of Mathematics, Williams College}
\email{\href{mailto:sjm1@williams.edu}{sjm1@williams.edu}}

\author{Vismay Sharan}
\address{Department of Mathematics, Yale University}
\email{\href{mailto:vismay.sharan@yale.edu}{vismay.sharan@yale.edu}}

\author{Marie-H\'el\`ene Tom\'e}
\address{Department of Mathematics, Duke University}
\email{\href{mailto:mariehelene.tome@duke.edu}{mariehelene.tome@duke.edu}}
\date{}

\numberwithin{equation}{section}

\begin{document}
\begingroup
\allowdisplaybreaks

\numberwithin{equation}{section}
\begin{abstract}
For a fixed elliptic curve $E$ without complex multiplication, $a_p \coloneqq p+1 - \#E(\mathbb{F}_p)$ is $O(\sqrt{p})$ and $a_p/2\sqrt{p}$ converges to a semicircular distribution. Michel proved that for a one-parameter family of elliptic curves $y^2 = x^3 + A(T)x + B(T)$ with $A(T), B(T) \in \mathbb{Z}[T]$ and non-constant $j$-invariant, the second moment of $a_p(t)$ is $p^2 + O(p^{{3}/{2}})$. The size and sign of the lower order terms has applications to the distribution of zeros near the central point of Hasse-Weil $L$-functions and the Birch and Swinnerton-Dyer conjecture. S. J. Miller conjectured that the highest order term of the lower order terms of the second moment that does not average to zero is on average negative. Previous work on the conjecture has been restricted to a small set of highly nongeneric families. We create a database and a framework to quickly and systematically investigate biases in the second moment of any one-parameter family. When looking at families which have so far been beyond current theory, we find several potential violations of the conjecture for $p \leq 250,000$ and discuss new conjectures motivated by the data.
\end{abstract}

\maketitle

\section{Introduction}
We assume the reader is familiar with elliptic curves. For detailed references, see, for example, \cite{silverman2, silverman1}. Let $\mcE \rightarrow \mathbb{P}^1$ be a non-split elliptic surface over $\mathbb{Q}$ with Weierstrass equation
    \begin{align}
        \mcE\colon\ y^2\ =\ x^3 + A(T)x + B(T)
    \end{align}
with $4A(T)^3 + 27 B(T)^2 \neq 0.$ We can take $A(T), B(T) \in \mathbb{Z}[T]$ and its $j$-invariant is given by
    \begin{align} 
        j(T) \ \coloneqq \ 1728 \frac{4A(T)^3}{4A(T)^3 + 27B(T)^2}.
    \end{align}
If $0 \leq \max\{3\mathrm{deg}(A(T)), 2\mathrm{deg}(B(T))\} < 12$, then $\mcE$ is a rational surface. In addition, almost all specializations $T=t \in \mathbb{Z}$ result in an elliptic curve $E_t$,
    \begin{align}
        E_t \ \coloneqq \ y^2\ =\ x^3 + A(t)x + B(t),
    \end{align}
and we say the set of all such $E_t$ forms a one-parameter family of elliptic curves of $\mathbb{Q}$. The rank of such a family is defined to be the minimum rank $r$ that appears infinitely often among the curves $E_t$ in the one-parameter family. 
The expected value of $\# E_t(\mathbb{F}_p)$ is $p+1$ and we write $a_t(p) \coloneqq p + 1 - \# E_t(\mathbb{F}_p)$ for the trace of Frobenius of $E_t$ at $p$. By Hasse's theorem on elliptic curves \cite{Hasse} we have that $|a_t(p)| \leq 2\sqrt{p}$.

The $n$-th moment of the Frobenius trace of a one parameter family is defined for each prime $p$ as
    \begin{align*}
        \mcA_{n, \mcE}(p) \ \coloneqq \ \sum_{t = 0}^{p-1} a_t(p)^n.
    \end{align*}
Note that we do not normalize by $1/p$, so this sum is always an integer. The moments of a one-parameter family encode arithmetic properties of the family. For example, Rosen and Silverman \cite{Rosen} proved a conjecture of Nagao \cite{Nagao} relating the first moment to the rank of an elliptic surface: if Tate's conjecture holds for the elliptic surface $\mcE$ (e.g., when $\mathcal{E}$ is a rational surface; see \cite{Shioda}), then
    \begin{align}
        \lim_{X \rightarrow \infty} \frac{1}{X} \sum_{p \leq X} \frac{\mcA_{1,\mcE}(p)\log p}{p} \ = \ - \mathrm{rank}\mcE(\mathbb{Q}(T)).
    \end{align}
Therefore, the rank of such a surface is determined by its first moments, and is directly related to negative bias in the values of $a_t(p)$.

Motivated by the negative bias in the first moment and its arithmetic importance, it is natural to ask if biases also exist in higher moments and if these biases are also arithmetically significant. This is the impetus for our work: we build a large database of Frobenius traces for primes up to 250,000 so that we can then easily investigate any one-parameter family, by, for example, numerically computing higher moments and attempting to isolate their lower order biases. Given more computing power, the framework we establish can easily be used to generate a larger database, and investigate higher moments.

To make our investigations precise, we recall the following asymptotic second moment expansion.
\begin{thm}[Michel \cite{michel}]\label{michels}
    Let $\mcE$ be an elliptic surface with nonconstant $j$-invariant. Then the second moment is of the form
        \begin{align}
            \mcA_{2, \mcE}(p) \ = \ p^2 + O\left(p^{3/2}\right),
        \end{align}
    with lower order terms of size $p^{3/2},\ p,\ p^{1/2}$ and $1$, respectively, where each has a cohomological interpretation.
\end{thm}
In his Ph.D. thesis \cite{Miller}, S. J. Miller computed a closed form expression for the second moment of several carefully chosen families. Based on these computations, and motivated by biases in the first moment expansion, Miller formulated the Bias Conjecture for the second moment: the largest lower order term in the second moment expansion which does not average to 0 is on average negative. As an application, negative biases in the second moment expansion are related to low-lying zeros and thus the average rank of $\mcE$ over $\mathbb{Q}(T)$, as this latter quantity can be bounded using the density of low-lying zeros. See \cite{SMALL2023, Miller05} for a detailed discussion to the excess rank problem.

Further results confirming the Bias Conjecture have been obtained in \cite{Asada, KN1, KN2, MMRW16, Miller, Miller04, Miller05}. For instance, see \cite{SMALL2023} for a table of families with closed form first and second moment expressions as obtained by \cite{Asada, Miller, Miller05}. The majority of these calculations rely on linear and quadratic Legendre symbol identites which we briefly touch on later; the takeaway is that these calculations are only tractable for lower degree $A(T),B(T)$, although this restriction often is still not enough. However, some results have also been obtained for specially chosen cubics and quartics. For example \cite{SMALL2023} recently proved that for primes $p \equiv 2\bmod 3$, the family $\mcE: y^2 = x^3 + x + T^3$ has second moment $p^2 + p$ by making use of an automorphism to lower the degree of $B(T)$. For primes $p \equiv 1\bmod 3$, however, they were not able to find a closed form expression and the data here appears random. Of note, the overall data appears generally slightly negative; we touch more on this later.

Nonetheless, this demonstrates progress towards disproving the Bias Conjecture, as a positive density of primes have positive bias.

As such, by looking at families  beyond the scope of Legendre symbol identities, we already begin to see interesting behavior. The salient point is that all of the families for which we have a closed form expression are specially chosen to be parameterized by low degree polynomials $A(T)$ and $B(T)$. Hence, we do not expect them to necessarily reflect the behavior of generic families. To this end, we computationally investigate one-parameter families defined by higher degree polynomials $A(T)$ and $B(T)$ where Legendre symbol calculations are intractable and to which the results of \cite{KN1, KN2} on cubic pencils do not apply. Our numerical work seeks to evaluate whether we expect the Bias Conjecture to hold by investigating a more generic swath of families. To isolate the lower order bias, we compute the following normalized second moment
    \begin{align}
        \mcB_{2, \mcE}(p) \ = \ \frac{\mcA_{2, \mcE}(p) - p^2}{p^{3/2}}
    \end{align}
and take a running average over primes $p$. Note that if the largest lower order term which does not average to $0$ is of size $p$ or $p^{1/2}$, then isolating lower order biases is difficult due to noise from the $p^{3/2}$ term drowning out fluctuations on the order of $p^{-1/2}$. Thus studying the Bias Conjecture numerically is a challenging problem. While our investigations turn up interesting families for further study, these techniques cannot prove or disprove the Bias Conjecture and are most useful for determining where to best allocate time and effort for further theoretical study and motivating the development of conjectures. Indeed, from a theoretical point of view, higher moments are signficantly more complicated to study, yet from a computational point of view, using our database one can compute them equally as easily. Hence, our work can be used to inform hypotheses about lower order biases in higher moment expansions.

We expect our database to have many applications to computing other quantities related to the traces of Frobenius of a one-parameter family and can be used to attempt to formulate analogous Bias Conjectures for the higher moments $n \geq 3$. In Section \ref{sec:database}, we discuss the techniques used to optimize computation and storage of the $a_p$ values to generate a database storing data for a large number of primes. Then, in Section \ref{sec:fams}, we illustrate the utility of the database by performing an extensive family search for polynomials $A(T)$ and $B(T)$ of degree at most $5$ and some of degree at most $10$. We investigate whether these families have potential positive bias using two statistics, an unweighted running average of the normalized second moment and a log-weighted average of the normalized second moment. Our search found some candidates with potential positive bias for further investigation which we also detail in this section. We also numerically explored the distribution of the normalized second moment and we discuss our findings and a conjecture about the variance of these distributions which our numerical work generated in Section \ref{sec:Conjectures}. Finally, in Section \ref{sec:future work}, we propose additional avenues for future work.

\section{Creating a database of \texorpdfstring{$a_p$}{} values}\label{sec:database}
There exists an explicit formula to compute $a_t(p)$ by calculating the number of solutions to $E_t\bmod p$ using Legendre symbols. For $a\in \mathbb{F}_p,$ the Legendre symbol is defined as
\begin{equation}
    \left( \frac ap \right)\ =\ \begin{cases}
        0 & a=0, \\
        1 & a=n^2 \text{ for some $n$},\\
        -1 & \text{otherwise.}
    \end{cases}
\end{equation}
Since $a_t(p) = p+1-\# E_t(\mathbb{F}_p),$ if $E_t: y^2 = x^3+A(t)x+B(t)$, then {for $p$ odd} we have that 
\begin{align}
    -a_t(p)\ =\ \sum_{x = 0}^{p-1} \left(\frac{x^3+A(t)x+B(t)}{p}\right).
\end{align}
Hence the second moment can be computed as
\begin{align}
    \mcA_{2, \mcE}(p) \ = \ \sum_{t=0}^{p-1}\sum_{x = 0}^{p-1}\sum_{w=0}^{p-1} \left(\frac{x^3+A(t)x+B(t)}{p}\right)\left(\frac{w^3+A(t)w+B(t)}{p}\right).\label{secondmomentformula}
\end{align}
For the families for which Miller was able to obtain a closed form expression, switching the order of summation yielded a quadratic or linear Legendre sum which can be evaluated using the following lemma. Then, for the carefully chosen families he considered, the remaining sums were easy to evaluate.
\begin{lemma}[Linear and Quadratic Legendre Sums]
    Let $p$ be an odd prime. If $a \not\equiv 0\bmod p$, then
        \begin{align}
            \sum_{t=0}^{p-1} \left(\frac{at + b}{p}\right)\ =\ 0,
        \end{align}
    and 
        \begin{align}
            \sum_{t=0}^{p-1} \left(\frac{at^2 + bt + c}{p}\right)\ =\ \begin{cases}
                (p-1)\left(\frac{a}{p}\right) &\text{\emph{if }} p \mid (b^2 - 4ac), \\
                -\left(\frac{a}{p}\right) &\text{\emph{otherwise.}}
            \end{cases}
        \end{align}
\end{lemma}

From \eqref{secondmomentformula}, we compute a database of 
$a_p$ for elliptic curves $E_{a,b}: y^3 = x^3 + ax + b$ for all values of $a,b \in \mathbb{F}_p$ for primes $p$ up to some prime $P$. This
requires both efficient algorithms to cut down computation time and efficient storage to minimize the amount of data that needs to be stored for each prime $p$. 
Naively, we need to compute $p^2$ values of $a_p$, however, in our database, 
we need only at most $4p+6$ values: namely,
    \begin{align}
        \# \mathrm{ \ of \ } a_p \mathrm{ \ stored \ } \ = \ \begin{cases} 4p & p\equiv 1\mathrm{ \ mod \ } 4 \\
    2p & p\equiv 3\mathrm{ \ mod \ } 4 \end{cases} \ + \ \begin{cases} 6 & p\equiv 1\mathrm{ \ mod \ } 3 \\
    0 & p\equiv 0,2\mathrm{ \ mod \ } 3. \end{cases}
    \end{align}
Using an additional optimization which we did not implement but which we detail below, one can reduce this to $2p+\theta(1)$ many $a_p$ for all primes $p$. In building the database used in this paper, we take $P$ to be the largest prime smaller than 250,000. However, our computation, storage, and look-up methods and can be extended to build a database for larger values of $P$ and extract values of $a_p$ for any $E_{a,b}$ corresponding to a pair of residues mod $p$ for any prime $p \leq P$ (up to limits in accurately storing large integers). The \verb|C++| code for computing and storing the values of $a_p$ for each $p$ is available at \cite{github} and is contained in the file \verb|quarticclasses.cpp|. Note that we store and compute $-a_p$ rather than $a_p$.

\subsection{Reduction to a smaller set of \texorpdfstring{$a_p$}{}'s}
We take advantage of two automorphisms and Legendre symbol identities to reduce the number of $a_p$ we need to compute and store. Since we are working over a field with characteristic zero, we write our elliptic curve with Weierstrass form $E: y^2 = x^3 + ax + b$ for some $a, b \in \mathbb{Q}$. By multiplying by an appropriate denominator, we take $a,b\in\mathbb{Z}.$ 

We reduce the number of $a_p$'s necessary by working with only certain residue classes of $a$ through a standard Weierstrass substitution. Specifically, for $\ell \in \ff_p^\times$, we let $y = \ell^3 \tilde{y}$ and $x = \ell^2 \tilde{x}$ and we obtain
    \begin{align}
        \ell^6 \tilde{y}^2 \ = \ \ell^6\tilde{x}^3 + \ell^2 a \tilde{x} + b,
    \end{align}
so that dividing by $\ell^6$ we arrive at the elliptic curve $\tilde{E}: y^2 = x^3 + a\ell^{-4} x + b\ell^{-6}$. The elliptic curves $E$ and $\tilde{E}$ are automorphic to one another and thus have the same number of solutions. Now, if we compute the $a_p$ for each quartic (i.e., fourth power) residue class, we can now compute the $a_p$ values for any elliptic curve over $\ff_p$, i.e., it suffices to compute the values of $a_p$ for one $a$ in each quartic residue class. There are five quartic residue classes if $p \equiv 1\bmod 4$ and three quartic residue classes if $p \equiv 3\bmod 4$. In each case, one quartic residue class corresponds to $a = 0.$ In the case $a = 0$, we further reduce the number of $a_p$ that need to be computed and stored. 

When $a=0$ and $p \not\equiv 1\bmod 3$, then $x \mapsto x^3$ is a bijection. Hence
    \begin{align}
        \sum_{x = 0}^{p-1} \left(\frac{x^3 + b}{p}\right) \ = \ \sum_{x = 0}^{p-1} \left(\frac{x + b}{p}\right) \ = \ 0.
    \end{align}
When $a=0$ and $p \equiv 1 \bmod 3$, it suffices to store one value of $a_p$ for each sextic residue class of $b\bmod p$, of which there are six. Hence we have reduced to the following:
    \begin{align*}
        \# a_p \mathrm{ \ we \ need \ to \ compute \ and \ store} \ = \ \begin{cases}
            4p &\text{if \ } p \equiv 1 \bmod 4 \\
            2p &\text{if \ } p \equiv 3 \bmod 4
        \end{cases} \ + \ \begin{cases}
            6 &\text{if \ } p \equiv 1 \bmod 3\\
            0 &\text{otherwise.}
        \end{cases}
    \end{align*}

The following is not implemented in our construction of the database but can be used to further reduce the number of $a_p$ values that need to be computed and stored. Suppose $-1$ is a quadratic residue in $\mathbb{F}_p$ (i.e., $p\equiv 1 \bmod 4$) and let $i \in \mathbb{F}_p$ be such that $i^2 \equiv -1\bmod p$. Then, make the substitutions $y = i\tilde{y}$ and $x = -\tilde{x}$ so that
    \begin{align}
        y^2 \ &= \ x^3 + ax + b\nonumber \\
        (i\tilde{y})^2 \ &= \ -\tilde x^3 - a\tilde x + b\nonumber \\
        \tilde{y}^2 \ &= \ \tilde x^3 + a\tilde x - b.
    \end{align}
Thus, when $p \equiv 1\bmod 4$, the curve $y^2 = x^3 + ax + b$ and the curve $y^2 = x^3 + ax - b$ are automorphic so it suffices to compute the values of $a_p$ for $0 \leq b < p/2$. This means we now have that
    \begin{align*}
        \# a_p \mathrm{ \ we \ need \ to \ compute \ and \ store} \ = \ \begin{cases}
            2p + 2 &\text{if \ } p \equiv 1 \bmod 4 \\
            2p &\text{if \ } p \equiv 3 \bmod 4
        \end{cases} \ + \ \begin{cases}
            6 &\text{if \ } p \equiv 1 \bmod 3\\
            0 &\text{otherwise.}
        \end{cases}
    \end{align*}

To retrieve the value of $a_t(p)$ corresponding to $E_t\colon\  y^2 = x^3 + A(t)x + B(t)$, we first find the quartic residue class of $A$. If $A \not\equiv 0\bmod p$, then for $a = 1,2,\ldots,$ we compute $A(t) a^{-1}\bmod p$ and check whether it is a fourth power mod $p$. Once we have found some $a$ with $A(t) a^{-1} \equiv \ell^4\bmod p$, we compute the corresponding value of $B(t)$ as $B(t) \ell^{-6}\bmod p$ and retrieve the $a_p$ value corresponding to $(a,b)$ (i.e., the automorphic curve $y^2 = x^3 + ax + b = x^3 + A(t)\ell^{-4}x + B(t)\ell^{-6}$). Recall that if $p \not\equiv 1\bmod 3$ and $A(t)$ is zero mod $p$, then $a_p = 0$ for all values of $B(t)$. Likewise, for all primes $p$, if $A(t)$ and $B(t)$ are both zero mod $p$ then $a_p = 0$. If $p \equiv 1\bmod 3$ and $B(t)\bmod p$ is nonzero, we find the sixth power residue class of $B(t)$ by computing $B(t) b^{-1}\bmod p$ for each value of $b$ stored until this is a sixth power residue. Once we find $b$ such that $B(t) b^{-1} \equiv \ell^6\bmod p$, we retrieve the value of $a_p$ for the curve $y^2 = x^3 + b$.

The only non-trivial arithmetic computation which is necessary is computing square roots of $x\bmod p$ (given that $x$ is a square). When $p\equiv 3 \bmod 4$ we have 
\begin{align}
    (x^{\frac{p+1}{4}})^2\ =\ x^{\frac{p+1}{2}}\ =\ x^{\frac{p-1}{2}}x\ =\ x,
\end{align}
where $x^{\frac{p-1}{2}} = 1$, since $x$ is a square. Unfortunately, when $p\equiv 1 \bmod 4$ we must apply a slightly more complicated idea: Cipolla's algorithm (see, for example, \cite{cipolla} for the development of the algorithm). We first find some $a\in \zz$ such that $a^2-x$ is not a square $\bmod\ p.$ It is unknown if this can be done deterministically, however this can easily be done probabilistically extremely quickly. Since $a^2-x$ is not a square, we conclude that 
\begin{align}
    \ff_{p^2}\ \cong\ \ff_p(\sqrt{a^2-x}).
\end{align}
Now, we have
\begin{align}
    \left(a+\sqrt{a^2-x}\right)^p\ =\ a^p + \left(\sqrt{a^2-x}\right)^p\ =\ a + \sqrt{a^2-x}\left(a^2-x\right)^{\frac{p-1}{2}}\ =\ a-\sqrt{a^2-x}.
\end{align}
Thus, $(a+\sqrt{a^2-x})^{p+1} = x,$ and thus $((a+\sqrt{a^2-x})^{\frac{p+1}2})^2 = x.$ Notably, since $x$ is a square in $\ff_p,$ we know that $(a+\sqrt{a^2-x})^{\frac{p+1}2} \in \ff_p.$ We note that this algorithm can be computed efficiently by repeated squaring. An observant reader will note that when $p$ is $3\bmod 4$, one can take $a=0$ and this algorithm then reduces to the previous case.

\section{Computationally Investigating the Bias Conjecture}\label{sec:fams}
While we expect our database to have many applications, we are initially motivated by calculating the second moment of one-parameter families of elliptic curves. Hence, we showcase how our database allows us to systematically computationally investigate the Bias Conjecture in one-parameter families for which obtaining a closed form of the second moment through Legendre symbols is intractable. The main, i.e., the most computationally expensive, inputs to calculating the second moment of a one-parameter family up to a value $P$ is creating a database of all of the possible values of $a_t(p)$ for all $p \leq P$. 

Since a one-parameter family is parameterized by $A(t)$ and $B(t)$, to compute the second moment, we need to look up $a_t(p)$, the $a_p$ value for each curve $E_t: y^2 = x^3 + A(t)x + B(t)$ resulting from specializing $T = t$, i.e., for the pair $(a,b) = (A(t), B(t))$, and then sum over $t$ as $t$ ranges over all possible residues mod $p$. To isolate and study the bias in the lower order terms, we compute the normalized second moment $\mcB_{2, \mcE}(p) = (\mcA_{2,\mcE}(p) - p^2)/p^{3/2}$ for each prime $p$. We take $P$ to be the largest prime smaller than $250,000$ and we compute the running average 
    \begin{align}
        \frac{1}{\pi(P) - 1} \sum_{2 < p \leq P} \mcB_{2, \mcE}(p),
    \end{align}
where $\pi(P)$ is the prime-counting function. We also introduce the log-weighted running average of the normalized second moment:
    \begin{align} 
        \frac{1}{N_w(P)}\sum_{2 < p \leq P} \mcB_{2, \mcE}(p)\log p,
    \end{align}
where $N_w(P) \coloneqq \sum_{2 < p \leq P} \log p\sim P$. This is desirable since as can be seen in the graphs that follow, the running average for smaller primes is not representative of the long-term behavior.

\begin{remark}
    Since $N_w(P)$ and $\pi(P)$ both go to infinity, we can safely ignore any finite number of primes in the limit. Indeed, we may want to ignore the finite number of additional primes that divide the $j$-invariant, $A(T)$, or $B(T)$ or perhaps where some ramification occurs. Further study is required to deduce what is optimal. 
\end{remark}

\subsection{Potential positive bias families}
We conducted an exhaustive family search for a potential positive bias family by looking at one parameter families defined by all possible combinations of polynomials $A(T), B(T)$ of degree $\leq 5$ with coefficients in $\{0, 1\}$. We computed the second moment for each of the resulting families for primes up to $1,000$ and noted those families whose running average of the normalized second moment was positive more than $95\%$ of the time.\footnote{Future work could focus on developing better measures of whether a family should be suspected of having positive bias and hence investigated further. For example, excluding primes $p \leq P_{\mathrm{min}}$.} We then computed the second moment for those families which passed this initial filtering for primes up to $250,000$ and calculated the graphs of the running averages (see Figures \ref{fig:fam1}, \ref{fig:fam2}, \ref{fig:fam3}, \ref{fig:fam4}, \ref{fig:fam5}, \ref{fig:fam6}). As a point of comparison, we computed the second moment of the family $y^2 = x^3 + x + T^3$ for primes up to $250,000$ (see Figure \ref{fig:SMALL23family}). 
\begin{figure}[ht]
    \centering
    \includegraphics[width=1.0\linewidth]{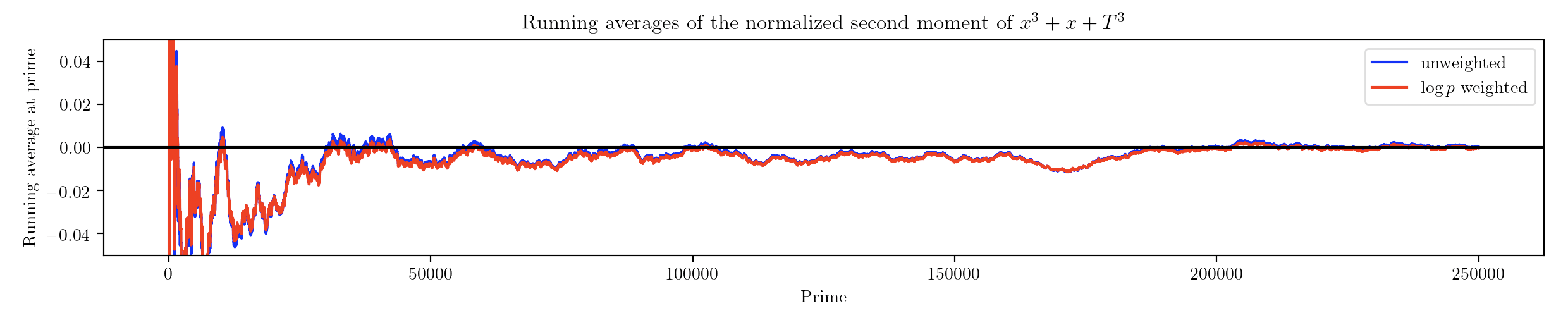}
    \caption{Running averages of the family investigated by \cite{SMALL2023}.}
    \label{fig:SMALL23family}
\end{figure}
Although the running average appears to stabilize after $200,000$, looking at the oscillations in the running average up to $200,000$ it is clear that one should still be concerned whether going out to $250,000$ is far enough to see the long range behavior of the running average. Hence, while it would also be beneficial to compute the $a_p$ values for a larger number of primes, this highlights the limitations of a computational approach. Our numerical investigations generate hypotheses and directions for future investigation but cannot prove or disprove the Bias Conjecture.

\begin{figure}[ht]
    \centering
    \includegraphics[width=.95\linewidth]{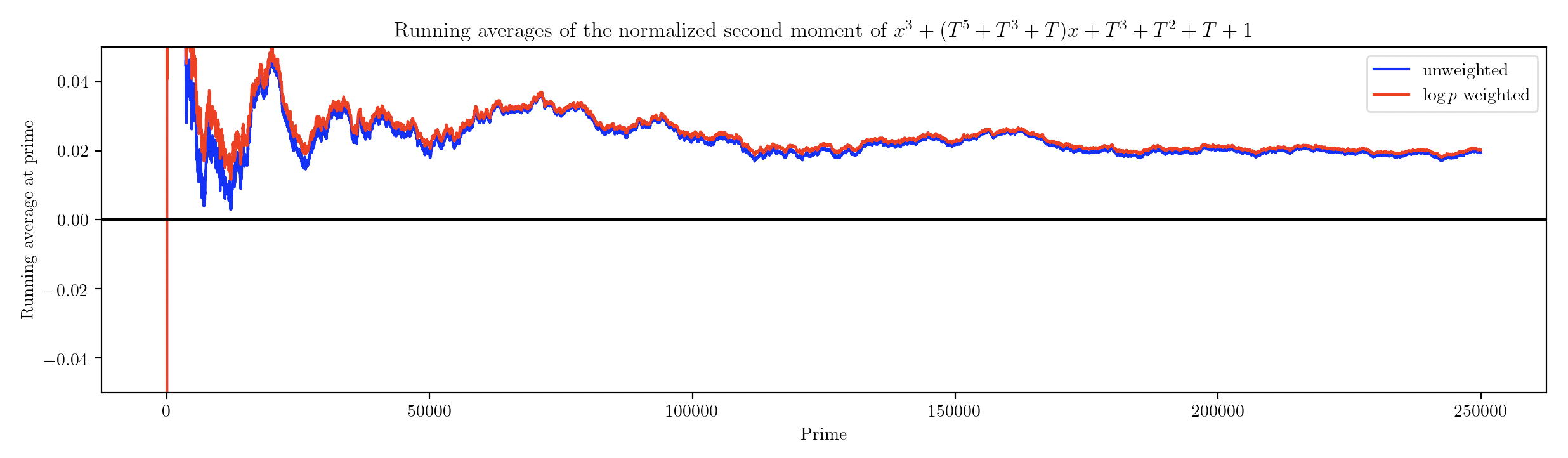}
    \caption{Running averages of $y^2\ =\ x^3 + (T^5 + T^3 + T)x + T^3 + T^2 + T + 1$.}
    \label{fig:fam1}
\end{figure}
\begin{figure}[ht]
    \centering
    \includegraphics[width=.95\linewidth]{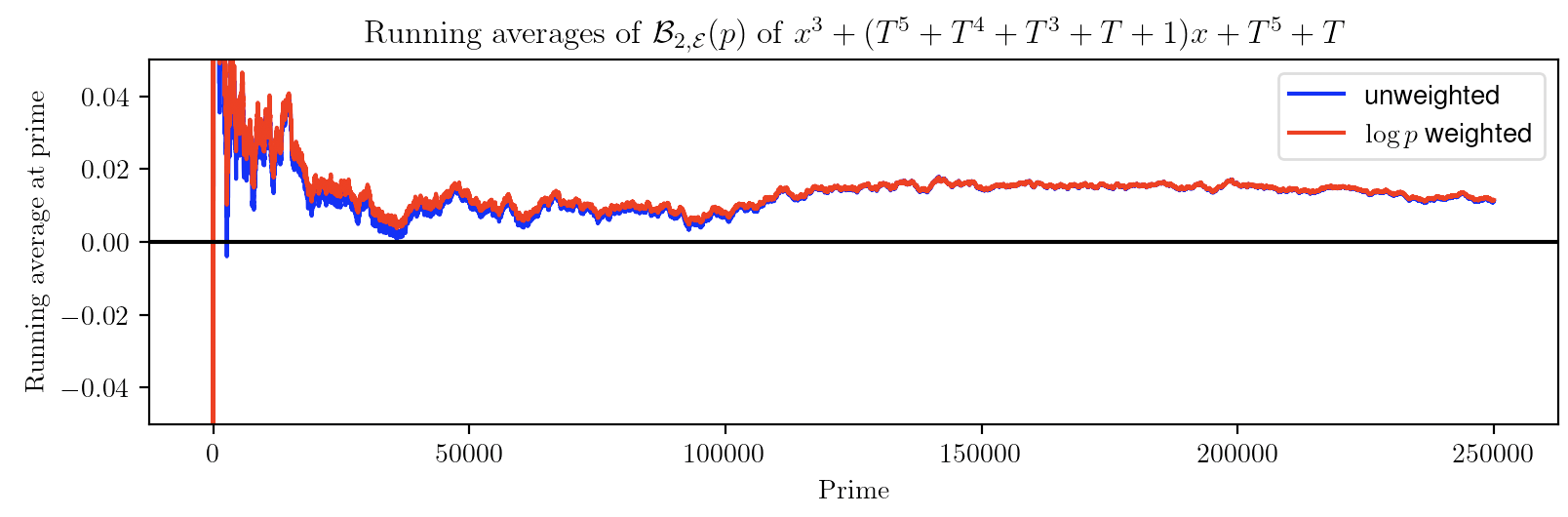}
    \caption{Running averages of $y^2\ =\ x^3 + (T^5 + T^4 + T^3 + T + 1)x + T^5 + T$.}
    \label{fig:fam2}
\end{figure}
\begin{figure}[ht]
    \centering
    \includegraphics[width=1.0\linewidth]{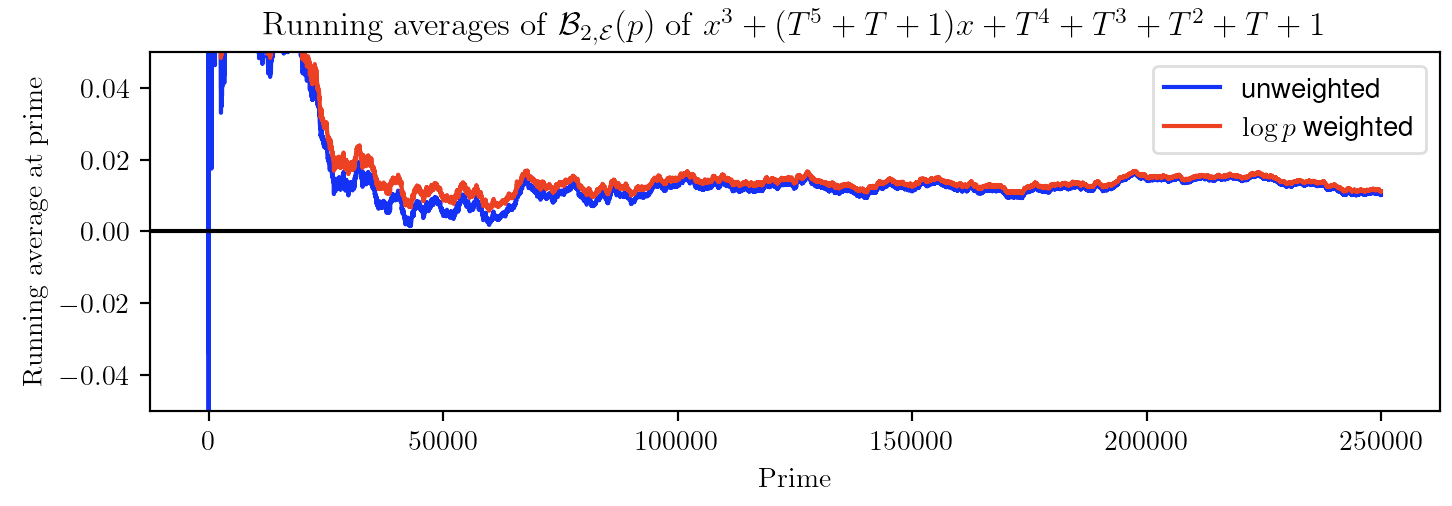}
    \caption{Running averages of $y^2\ =\ x^3 + (T^5 + T + 1)x + T^4 + T^3 + T^2 + T + 1$.}
    \label{fig:fam3}
\end{figure}
\begin{figure}[ht]
    \centering
    \includegraphics[width=1.\linewidth]{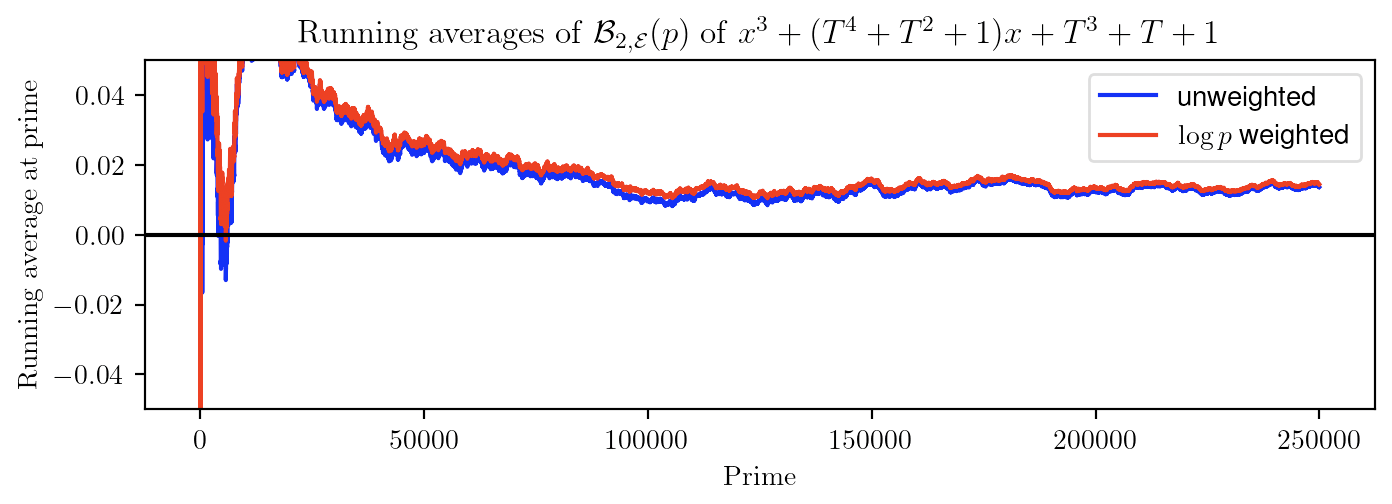}
    \caption{Running averages of $y^2\ =\ (T^4 + T^2 + 1)x + T^3 + T + 1$.}
    \label{fig:fam4}
\end{figure}
\begin{figure}[ht]
    \centering
    \includegraphics[width=1.\linewidth]{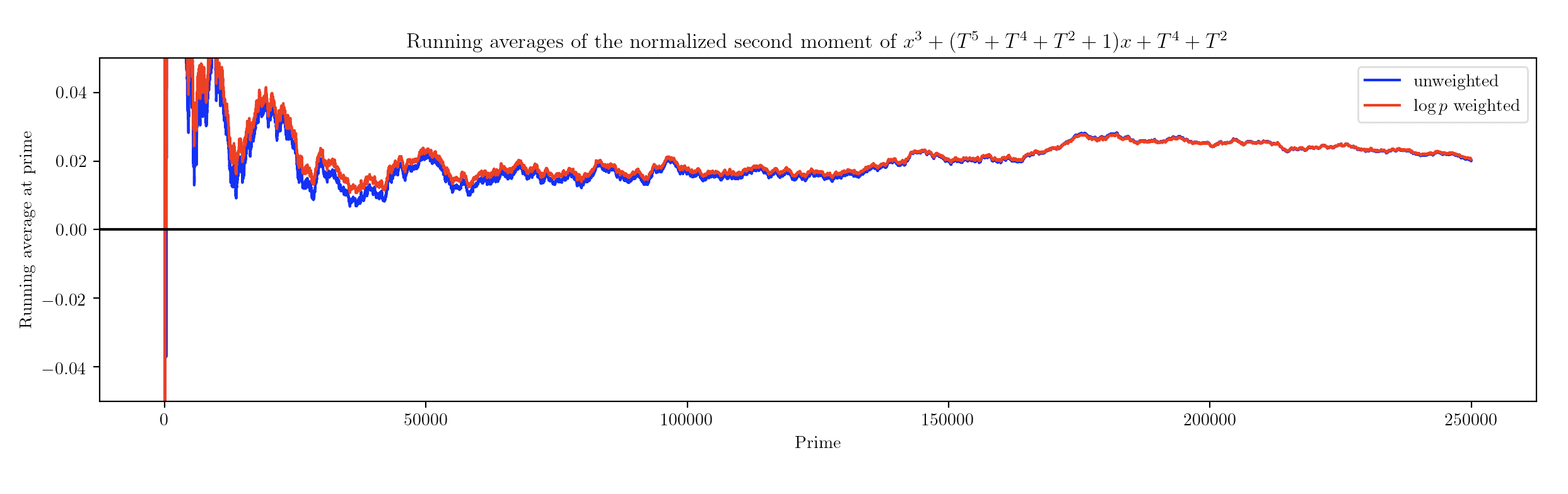}
    \caption{Running averages of $y^2\ =\ x^3 + (T^5 + T^4 + T^2 + 1)x + T^4 + T^2$.}
    \label{fig:fam5}
\end{figure}

\begin{figure}[ht]
    \centering
    \includegraphics[width=1.0\linewidth]{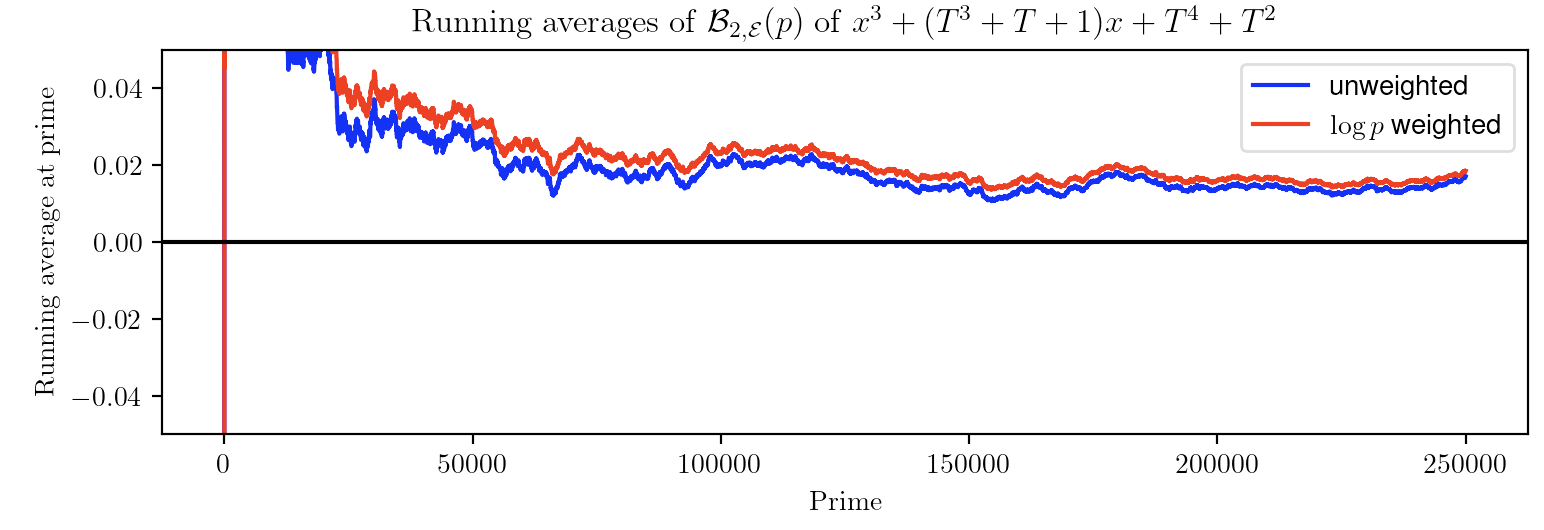}
    \caption{Running averages of $y^2\ =\ x^3 + (T^3 + T + 1)x + T^4 + T^2$.}
    \label{fig:fam6}
\end{figure}

It is instructive to compare the graphs of the running averages of the families our family search found to the family $y^2 = x^3 + x + T^3$ for which \cite{SMALL2023} proved that for half of the primes, the second moment of this family has positive bias. While the graph of the running averages for $y^2 = x^3 + x + T^3$ is negative most of the time but seems to tend to zero at $250,000$, the graphs in Figures \ref{fig:fam1}, \ref{fig:fam2}, \ref{fig:fam3}, \ref{fig:fam4}, \ref{fig:fam5}, \ref{fig:fam6} remain positive up to $250,000$ and seem to stabilize for primes beyond $100,000$. Therefore, with the above caveats in mind, we have reason to suspect these families may have positive bias.

\begin{remark}\label{arth prog}
When splitting up these familes based on the primes residue class mod $12$, we often see extremely different behavior in the resulting $4$ residue classes. Specifically, if $p$ is such that $A(T)$ and $B(T)$ split completely, the bias tends to be significantly higher. Because our polynomials' coefficients are either $0$ or $1,$ being $1\bmod{12}$ tends to result in splitting. 
\end{remark}
\newpage\newpage
\begin{remark}
    One may worry that positive graphs of the running averages the families we have chosen out of $2^{10}$ families explored are the result of random fluctuations rather than underlying structure in the one-parameter family which may yield positive bias in the lower order terms of the second moment. However, when factoring $A(T)$ and $B(T)$ which appear in the above families, $A(T)$ is either irreducible or has a factor of $T^2 - T + 1$, $T^2 + T + 1$, or $T^2 + 1$ and the same is true for $B(T)$. This leads us to suspect that the positive graphs of these families are due to more than just the random chance.   
\end{remark}

When looking at polynomials of degree at most $10$, we found the family $y^2 = x^3 + T^{10}x + T^8 + T^2$, which is our most promising candidate for having positive bias (see Figure \ref{fig:fammostpos}). At the largest prime smaller than $250,000$, the running average of the normalized second moment for this family is approximately $0.0287$ and the log-weighted running average is approximately $0.0283$, or a 4.2 $\sigma$ deviation. Note that the polynomial $B(T)$ has a factor of $T^2 + 1$. 

\begin{figure}[ht]
    \centering
    \includegraphics[width=1.0\linewidth]{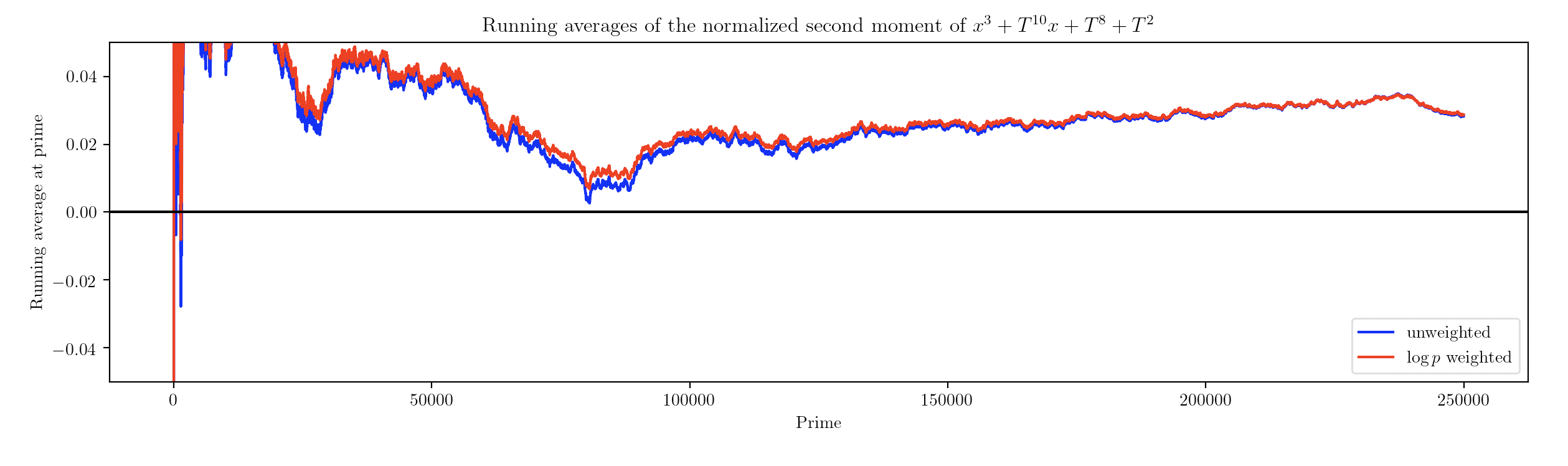}
    \caption{Running averages of $y^2\ =\ x^3 + T^{10}x + T^8 + T^2$}
    \label{fig:fammostpos}
\end{figure}

\section{Distribution of the Normalized Second Moment \texorpdfstring{$\mcB_{2, \mcE}(p)$}{}}\label{sec:Conjectures}
Motivated by the Sato-Tate Conjecture, we study the distribution of the normalized second moment of a one-parameter family. It is natural to ask what determines the variance of the distribution of the normalized second moment of a one-parameter family. Based on our numerical computations, we formulate the following conjecture.
\begin{conjecture}\label{var = int}
    The variance of the distribution of $\mcB_{2, \mcE}(p)$ always converges to a positive integer.
\end{conjecture}

Indeed there seems to be a deeper connection between the polynomials $A(T)$ and $B(T)$ and the conjectured integer the variance converges too. Indeed, the one parameter family given by $y^2 = x^3+A(T^n)x+B(T^n)$ seems to be an integral multiple of $y^2 = x^3+A(T)x+B(T)$ based on the prime factorization of $n,$ $A$ and $B.$ 

The following figures have binned the data between the maximum and the minimum over all primes less than $250,000$ and then split into 100 equal sized buckets. 

\begin{figure}[ht]
    \centering
    \begin{minipage}{0.47\textwidth}
        \centering
        \includegraphics[width=\textwidth]{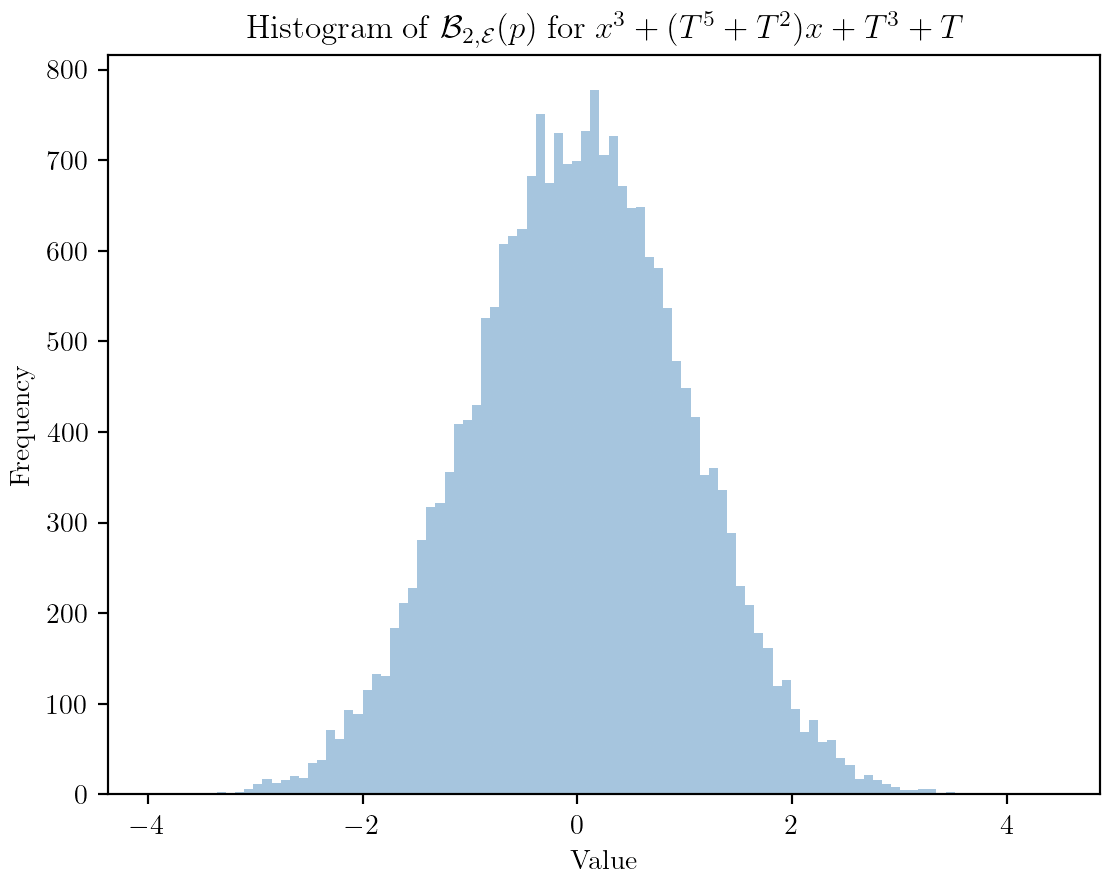}
        \caption{A Family with variance $1.015$}
        \label{fig:var1exa}
    \end{minipage}\hfill
    \begin{minipage}{0.46\textwidth}
        \centering
        \includegraphics[width=\textwidth]{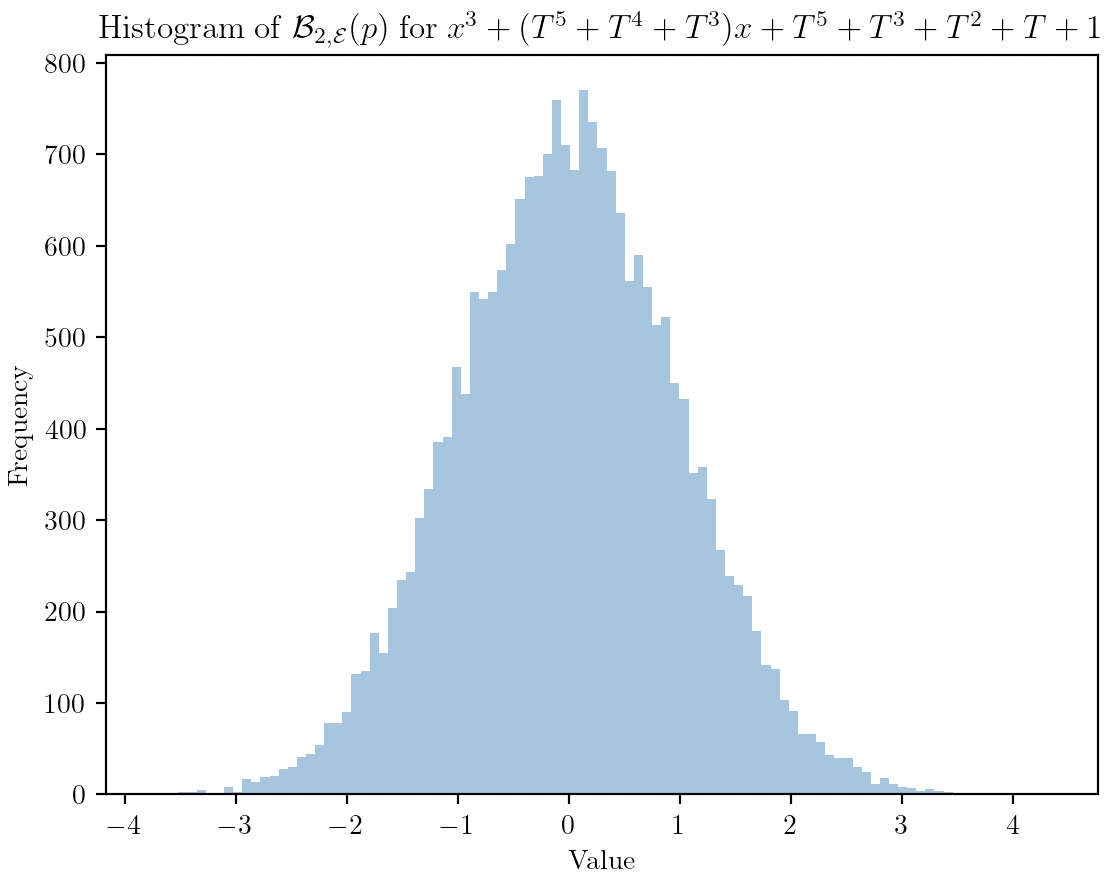}
        \caption{A Family with variance $1.005$.}
        \label{fig:var1exb}
    \end{minipage}
\end{figure} 

\begin{figure}[ht]
    \centering
    \begin{minipage}{0.45\textwidth}
        \centering
        \includegraphics[width=\linewidth]{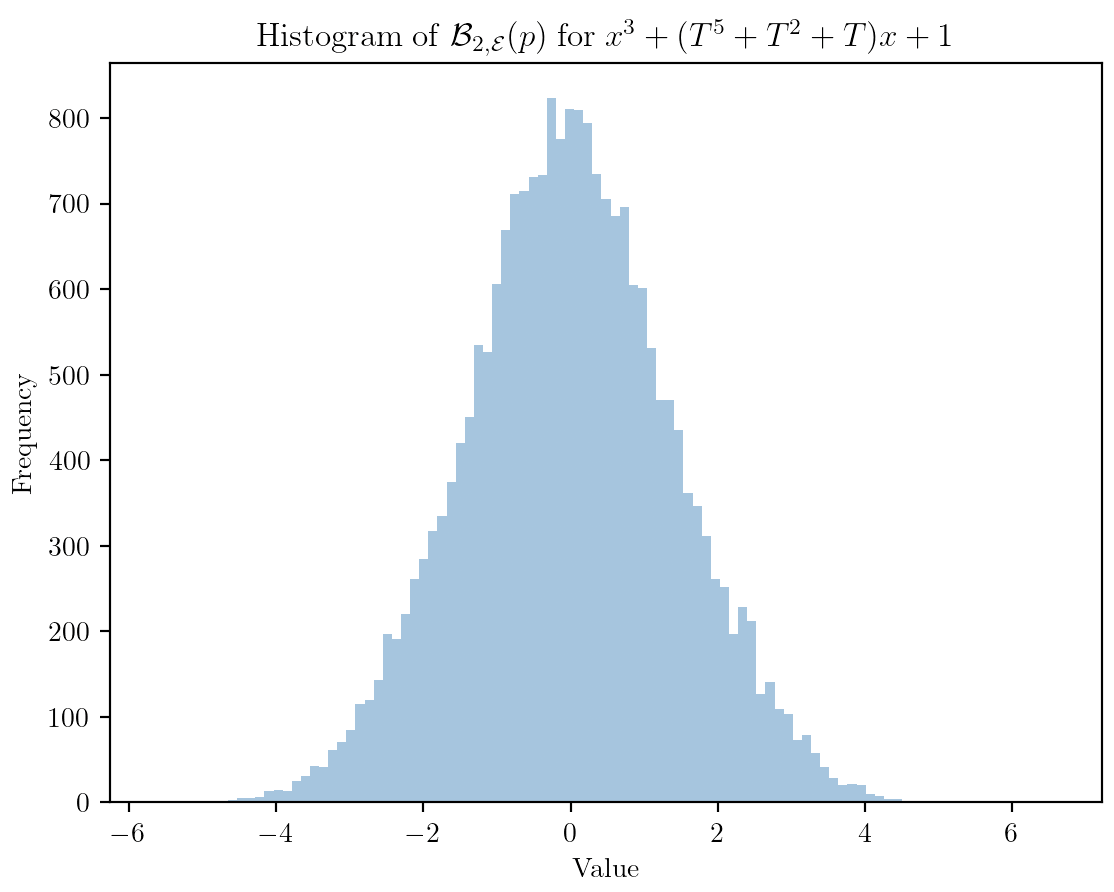}
        \caption{A family with variance $1.991$.}
        \label{fig:var2}
    \end{minipage}\hfill
    \begin{minipage}{0.45\textwidth}
        \centering
        \includegraphics[width=\linewidth]{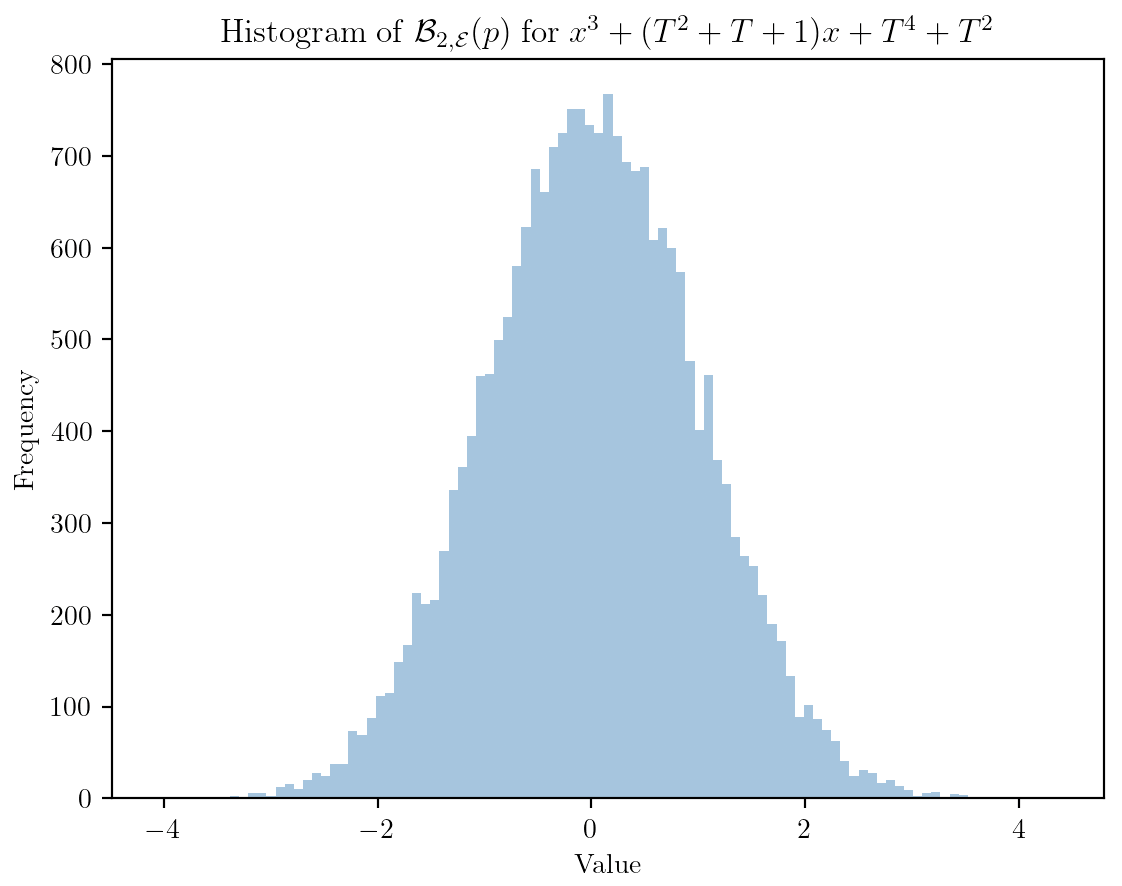}
        \caption{A family with variance $0.994$.}
        \label{fig:var1}
    \end{minipage}
\end{figure}

Our generic case seems to be variance converging to $1$ (see Figures \ref{fig:var1exa}, \ref{fig:var1exb},
and \ref{fig:var1}). In this generic case, the distribution of $\mcB_{2,\mcE}(p)$ seems to exhibit Gaussian-like behavior. We found one family $\mcE: y^2 = x^3 + (T^5 + T^2 + 1)x + 1$ whose variance seems to be converging to $2$ (see Figure \ref{fig:var2}). 

\begin{figure}[ht]
    \centering
    \includegraphics[width=\linewidth]{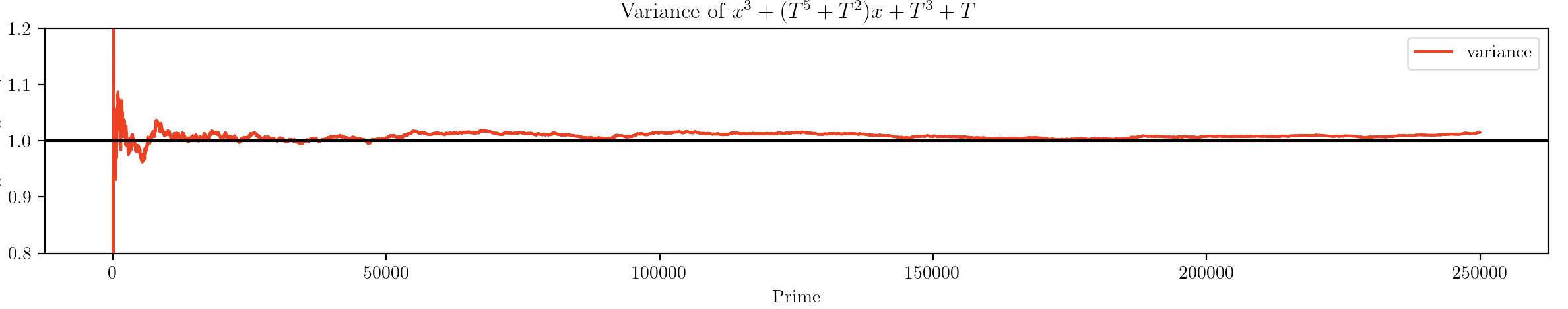}
    \caption{The graph of the variance of a family over primes $\leq p$.}
    \label{fig:vargraph}
\end{figure}

Additionally, when restricting the data to ignore the first $10,000$ primes, we net a family with the similar variance and distribution. Indeed, the only difference seems to be that the data is slightly more random, which is expected given that there are fewer data points.

\begin{figure}[ht]
    \centering
        \centering
        \includegraphics[width=.9 \linewidth]{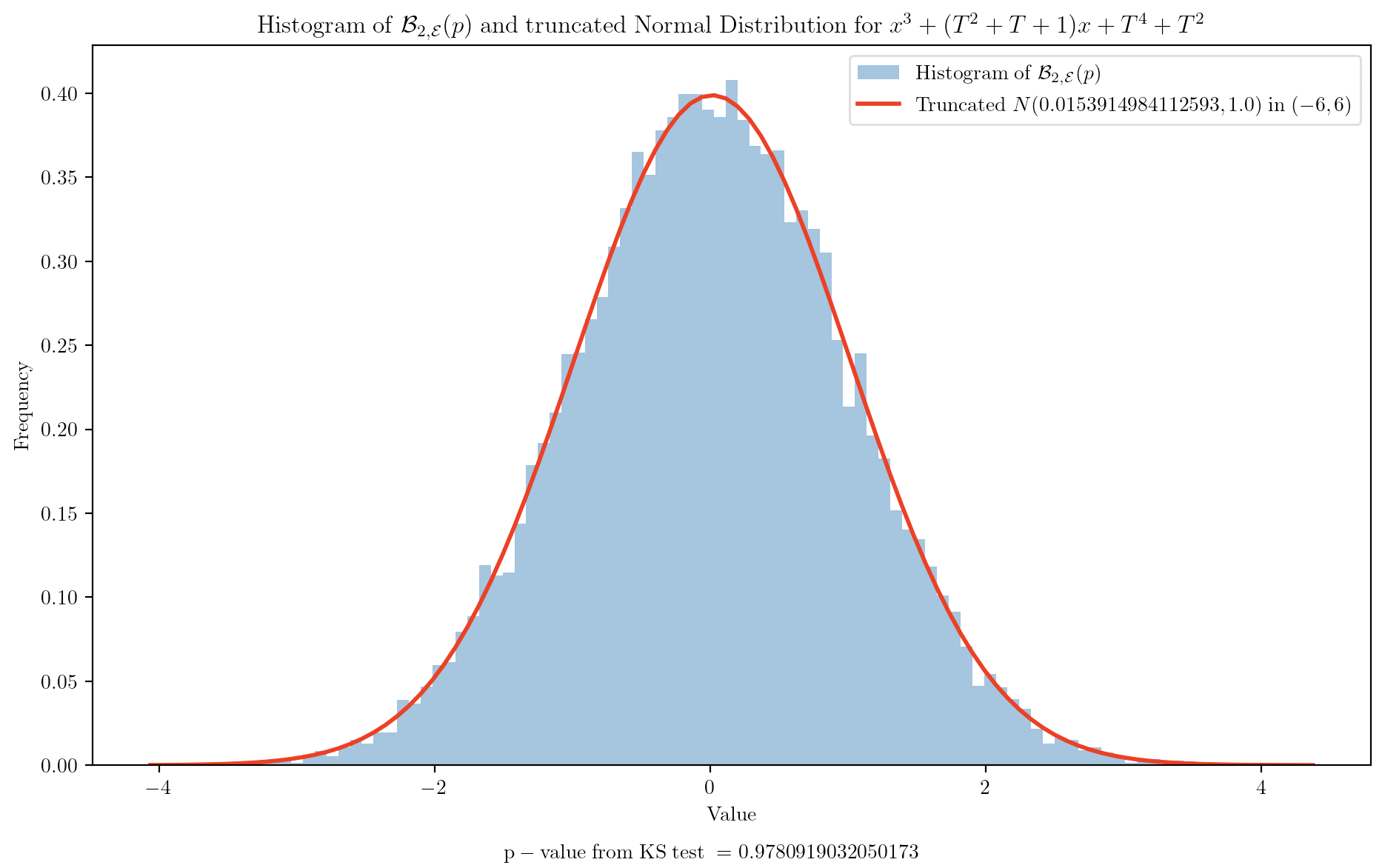}
        \vspace{-.3cm}
        \caption{Fitting a truncated normal to the variance $1$ family in Figure \ref{fig:var1}}
        \label{fig:var1fit}
\end{figure}

\begin{figure}[ht]
        \centering
        \includegraphics[width=.89\linewidth]{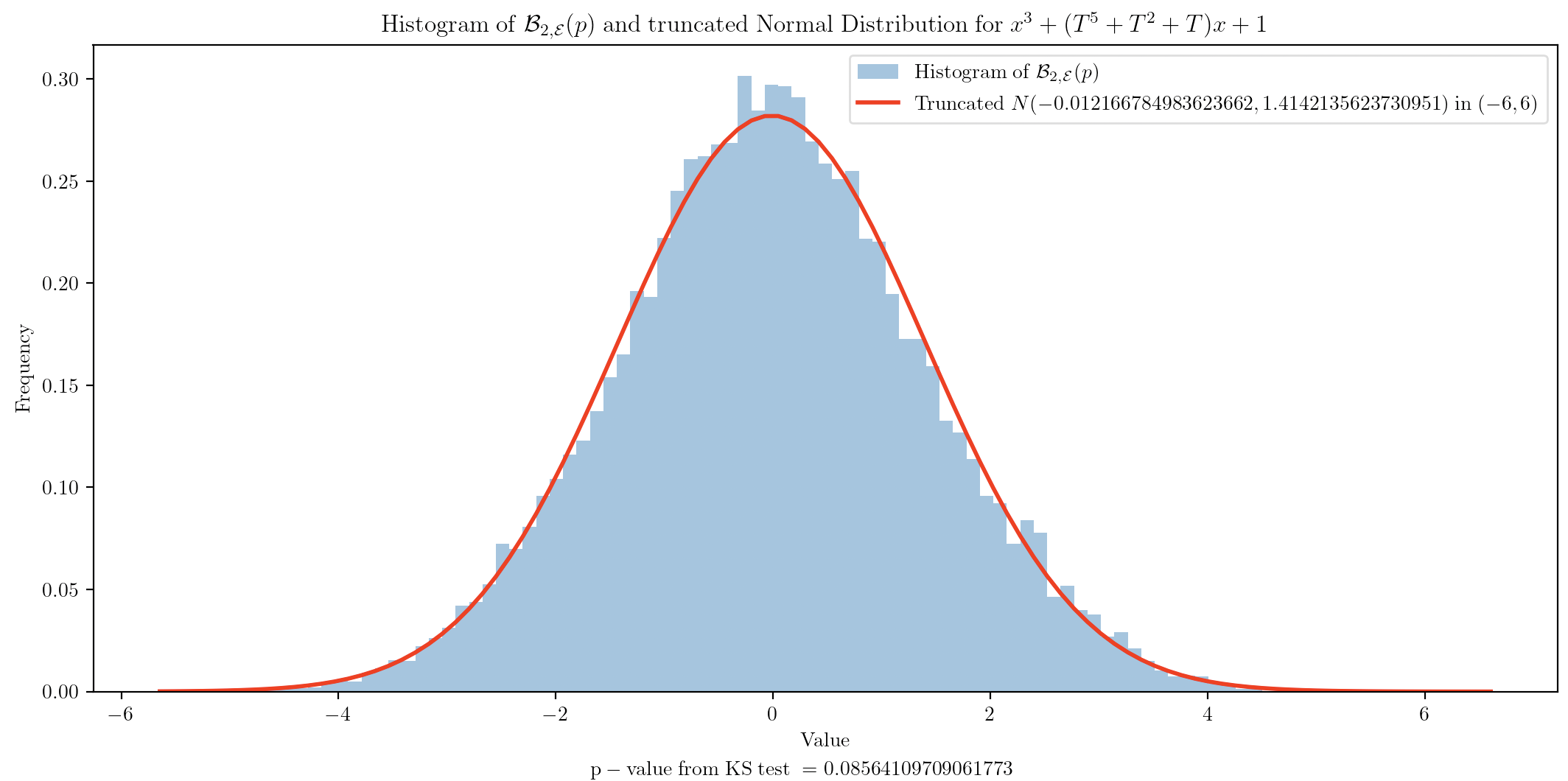}
        
        \caption{Fitting a truncated normal to the variance $2$ family.}
        \label{fig:var2fit}
\end{figure}

Thus, a natural next question is what distribution the normalized second moment of a one-parameter family converges to in our seemingly generic case of variance $1$ Gaussian-like behavior and what is our generic case, i.e., what restrictions can we impose on our one-parameter family to ensure we are in the generic case?

These distributions look like normal distributions, however they cannot be as they must be bounded due to Theorem \ref{michels}. This led us to investigate whether a truncated normal distribution with the same mean and variance as the family was a good fit to the data. When comparing the fit of a truncated normal distribution for a generic family to the variance $2$ family, the truncated normal is a better fit for the variance $1$ family while for the variance $2$ family, there is too much mass around zero (see Figures \ref{fig:var1fit} and \ref{fig:var2fit}).

Our numerics motivate the study of the distribution of the second moment of a one-parameter family. Further, the database allows for efficient computation of higher moments of one-parameter families. 

\section{Higher moments}
Calculating higher moments is difficult because it measures a more subtle distribution of $a_p$ values. However, numerically calculating higher moments is no different than the second moment. Indeed, one can calculate the second through tenth moments all at once nearly as quickly as just the second moment. 

To highlight some applications of our database, we calculate higher moments of a one-parameter family. Without the closed form for $a_p$ values we computed above, this method would not be possible. By Theorem 1 of \cite{birch},
\begin{align}
    \mcA_{2n, \mcE}(p)\ =\ C_np^{n+1} + O\left(p^{n+1/2}\right)
\end{align}
where $C_n = \frac{(2n)!}{n!(n+1)!}$ denotes the $n$th Catalan number. Like before, we investigate 
\begin{align}
    \mcB_{2n, \mcE}(p)\ \coloneqq\ \frac{(\mcA_{2n,\mcE}(p)/C_n) - p^{n+1}}{p^{n+1/2}}.
\end{align}
The (unweighted) running averages $\mcB_{2n,\mcE}$ for a generic family are shown below.

\begin{figure}[ht]
        \centering
        \includegraphics[width=\linewidth]{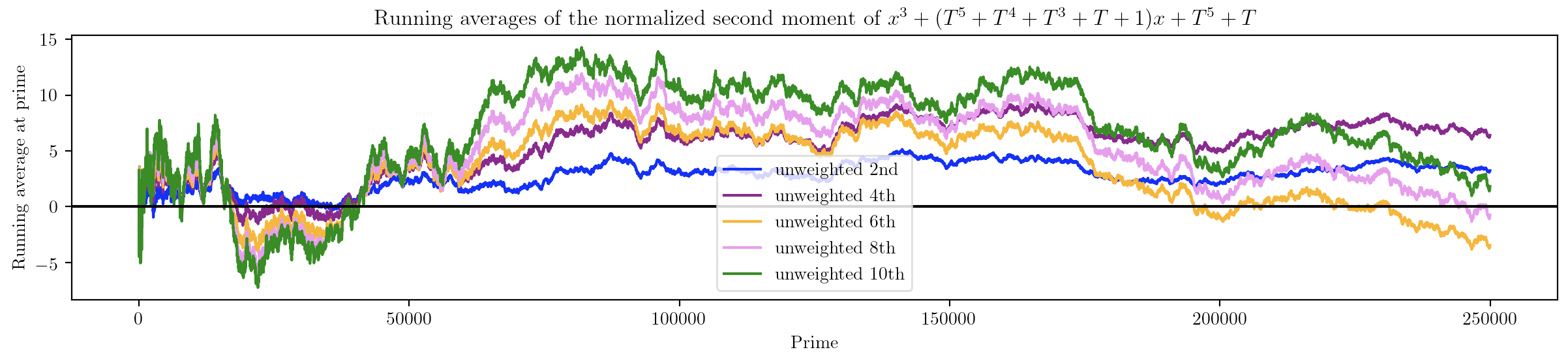}
        \vspace{-.3cm}
        \caption{Normalized 2nd through 10th moments}
        \label{fig:2ndthrough10th}
\end{figure}

\begin{figure}[ht]
        \centering
        \includegraphics[width=\linewidth]{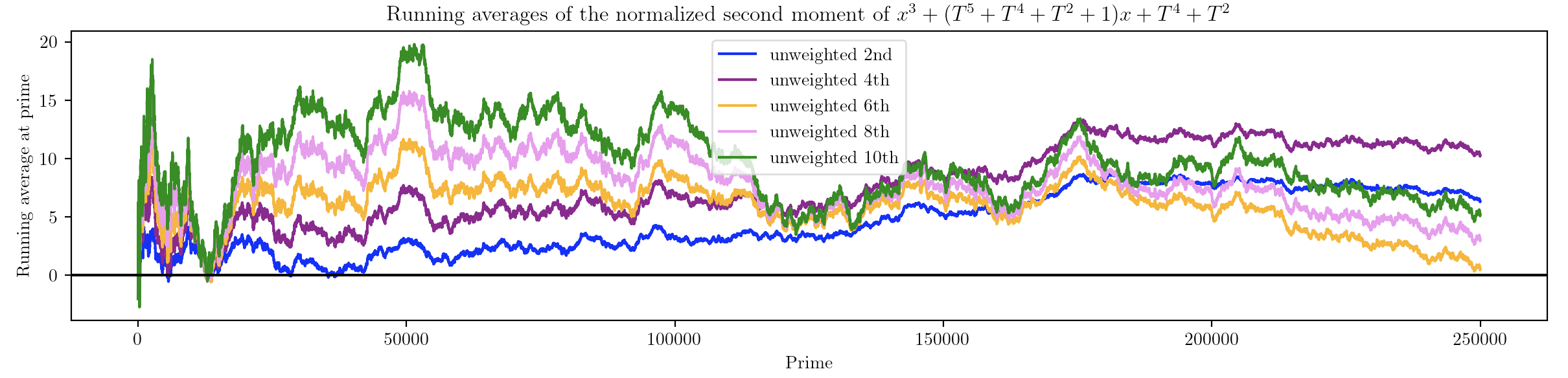}
        \vspace{-.3cm}
        \caption{Normalized 2nd through 10th moments}
        \label{fig:2ndthrough10th2}
\end{figure}

\begin{figure}[ht]
        \centering
        \includegraphics[width=\linewidth]{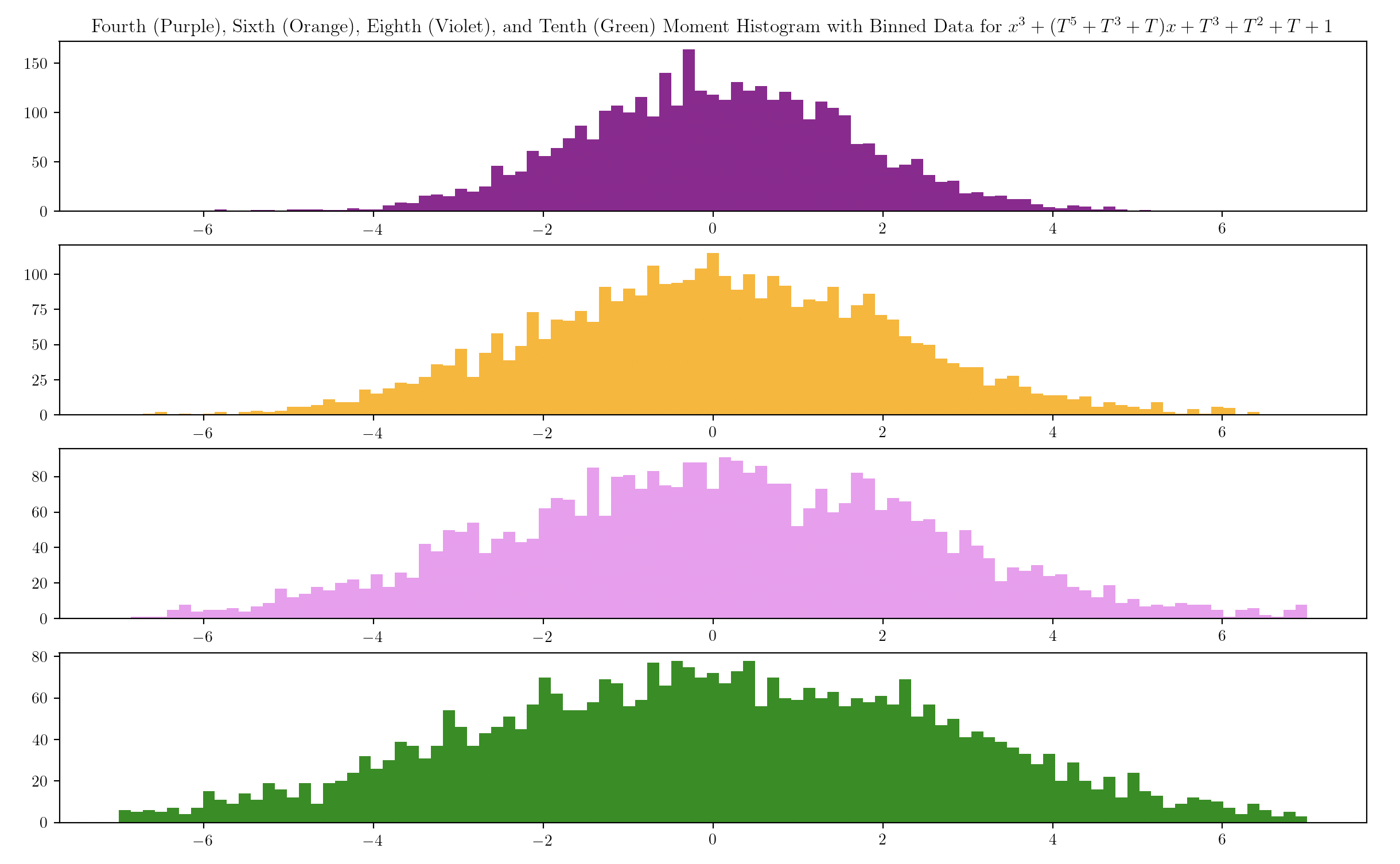}
        \caption{Normalized 4th through 10th moments, split into 100 buckets from -7 to 7}
        \label{fig:2ndthrough10thhist}
\end{figure}

We recall from above that in the second moment case, our variance converges to an integer, and the data is very similar to a normal distribution. We consider how these properties change for higher moments. Initially, the distribution of $\mcB_{4,\mcE}$ and $\mcB_{6,\mcE}$ look like normal distributions. There are a couple of notable differences. 

First, the distribution of $\mcB_{2n, \mcE}$ values seem to be spreading out as $n$ increases. This is equivalent to the variances increasing, and indeed the data seem to support this. One potential reason is that a small bias in $a_p$ values will result in higher deviations from the expectation in higher moments. 

Second, the variances no longer seem to be an integer. Perhaps this is because the correct way to normalize the $\mcA_{2n,\mcE}$ values is by 
\begin{align}
    \mcB'_{2n, \mcE}(p)\ \coloneqq\ \frac{\mcA_{2n,\mcE}(p) - C_np^{n+1}}{p^{n+1/2}}.
\end{align}
However, under this normalization, we do not have enough data to conclude that $\mcB'$ values either converge or do not converge to an integer.
\begin{remark}
    The difference between these two normalizations is that the second one multiplies the variance by $(C_n)^2.$
\end{remark}

Third, as $n$ increases, the distribution seems to stray further and further away from a smooth distribution. This suggests that convergence occurs at a slower rate for higher moments.

Finally, the higher moments do not seem to have a generic behavior of positive or negative bias: some families seem to have positive bias while others seem to have negative bias. Further investigation is warranted to understand how the polynomials $A(T)$ and $B(T)$ determine the value to which the higher moments converge.

\section{Future work}\label{sec:future work}
To compute moments up to a fixed number $X,$ we first need to compute approximately $2X^2/\log X$ many $a_p$ values, and then do a comparable amount of work to compute moments. 

However, computing $a_p$ values using the naive algorithm described above requires $O(p)$ computations. This algorithm's advantage is that due to summing over all of $\ff_p$ in some cases we can extract cancellation. However, if one simply wishes to compute $a_p$ values, one can employ Schoof's algorithm. See, for instance, \cite{schoofelkiesatkin} for a discussion of this algorithm which allows for computation of $a_p$ in $O(\log^6(p))$ times (assuming arithmetic operations are $O(1)$). Indeed, taking our $p=3,5,7,11,13$ one can compute $a_p$ values for primes less than $14,000,000$, which is approximately $50$ times as many primes as we computed. 

Computing significantly more $a_p$ values should provide computation evidence for Conjecture \ref{var = int}, allow for further testing of the bias conjecture, and probing higher moments, potentially even over certain residue classes mod $p.$ One can consider certain classes of primes, for instance, those which split, split completely, ramify etc. over a certain number field, to investigate how $A(T)$ and $B(T)$ may influence the behavior of the higher moments.

For potentially easier things to study, how does the variance, second moment, and distribution of second moments of $y^2 = x^3+A(T^n)x+B(T^n)$ relate to $y^2 = x^3+A(T)x+B(T)$? Additionally, what happens for primes where $A(T)$ and $B(T)$ both factor completely? Or perhaps, what happens when both $A(T)$ and $B(T)$ are irreducible. These questions dicussing behavior in these restricted environments provide fertile ground for explorations. We note that we have been splitting our data into $100$ buckets. One can explore how a different number of buckets changes how good the $p$-values from the KS test nets, and use this to work towards determining the distribution. 

Other potential areas of study include exploring higher moments, using techniques from algebraic geometry to study a threefold that holds second (or higher) moment information, finding the distribution that the second moment converges to, proving Conjecture \ref{var = int}, proving Remark \ref{arth prog}, and determining the integer that the second moment converges to based on the polynomials $A(T)$ and $B(T)$. 

\section*{Acknowledgements}
\addcontentsline{toc}{section}{Acknowledgements}
This research was completed at the Williams College SMALL REU Program and was supported by Williams College and the National Science Foundation (grant DMS2241623). The authors are grateful for the support of Duke University, Princeton University, and the University of Michigan. The authors would like to thank Adam Logan for helpful discussions.

\printbibliography
\end{document}